\documentclass[12pt]{article}
\usepackage[tbtags]{amsmath}
\usepackage{amsfonts,amssymb}
\usepackage{a4wide}
\usepackage{graphicx}

\begin{document}
\title{Small braids having a big Ultra Summit Set}
\author{Maxim Prasolov}

\maketitle

\section{Abstract.}

In [\ref{rigid}] authors asked: (open question 2) is the size of USS of a rigid pseudo-Anosov braid is bounded above by some polynomial in the number of strands and the braid length? We answer this question in the negative.

\section{Introduction.}

The conjugacy problem in the braid group is solved by F.Garside in 1969. The problem separates into two parts: to decide for two braids whether they are conjugate and to find by what braid they are conjugate. Advances in these problems` solutions are equal, so we do not differ them. In papers [\ref{rigid},\ref{periodic},\ref{bkl},\ref{elrifai},\ref{processing},\ref{franco},\ref{gebhardt}] the algorithm was improved but a polynomial algorithm is still unknown.

Garside`s algorithm calculates the Summit Set of a braid. We define it in Section 3. In short, the Summit Set is a finite canonical subset of the conjugacy class of the braid.

Later W.Thurston [\ref{processing}] improved the algorithm by introducing a left normal form of the braid which can be calculated in polynomial time in the number of strands and the braid length. So W.Thurston obtained a polynomial algorithm for the word problem in the braid group.

In 1994 E.El-Rifai and H.Morton [\ref{elrifai}] improved the algorithm by replacing a Summit Set by its subset Super Summit Set. And it was shown how one can find at least one element of the Super Summit Set in polynomial time.

And in 2003 J.Gonz\'alez-Meneses and N.Franco [\ref{franco}] shown how one can find the rest of the Super Summit Set in the time bounded above by a product of the size of the Super Summit Set and a polynomial in the number of strands and the braid length. However, the size of the Super Summit Set is not bounded above by a polynomial.

In 2003 V.Gebhardt [\ref{gebhardt}] replaced a Supper Summit Set by its subset Ultra Summit Set. The time to find at least one braid in the Ultra Summit Set is not known but one can find the rest of the Ultra Summit Set in the time bounded above by a product of the size of the Ultra Summit Set and a polynomial in the number of strands and the braid length. However, the size of Ultra Summit Set is not bounded above by a polynomial.

J.Birman, V.Gebhardt and J.Gonz\'alez-Meneses [\ref{rigid}] suggested reducing problem to so called rigid pseudo-Anosov braids. They asked: (open question 2) is the size of USS of a rigid pseudo-Anosov braid is bounded above by some polynomial in the number of strands and the braid length? In this paper we answer this question in the negative.

\smallskip

\noindent{\bf Theorem 1.} The braid $\alpha_n:=\sigma_1 \sigma_2^{-1} \sigma_3 \sigma_4^{-1}\dots \sigma_{n-1}^{(-1)^n}$ on $n$ strands is rigid and the size of its Ultra Summit Set is at least $2^{[(n-2)/2]}.$

\noindent{\bf Theorem 3.} The braid $\alpha_n$ is pseudo-Anosov for odd number $n\ge3.$

\smallskip

\noindent{\bf Note.} The statement of Theorem 3 holds true for arbitrary $n$. One can construct an invariant train-truck (for the notion of train truck see [\ref{bestvina},\ref{thurston}]) of the braid $\alpha_n.$ But for simplicity we present here a direct proof for odd $n$.

This counter-example was discovered experimentally by I.A.Dynnikov with help of the program of J.Gonz\'alez-Meneses which implements the last version of algorithm introduced in [\ref{gebhardt}]. According to the calculations the size of the Ultra Summit set of the braid $\alpha_n$ is equal to $$(3-(-1)^n)\cdot 3^{n-3},$$ if $n=3,4,\dots,11.$

\smallskip

J.Birman, K.H.Ko, S.J.Lee [\ref{bkl}] introduced a new presentation of the braid group with so-called {\it band generators}. Using it, they obtained a new solution to the word and conjugacy problem which retains most of desirable features of the solution we consider above (it uses the notion of left normal form and Ultra Summit Set) and at the same time makes possible certain computational improvements. However, Birman-Ko-Lee presentation modification of open question 2 in [\ref{rigid}] has the same answer.

\smallskip

\noindent{\bf Theorem 2} Ultra Summit Set in Birman-Ko-Lee presentation of the braid $\alpha_n$ contains at least $2^{(n{-}1)/2}$ rigid braids for odd $n$.

\smallskip

Due to computations of my program the size of USS of $\alpha_n$ is equal to $$\frac{3-(-1)^n}2\cdot n\cdot3^{n-3},$$ if $n=3,4,\dots,9.$

\section{Definitions.}

The following definitions refer to theorem 1.

\noindent{\bf Definition.} A braid with positive crossings is a {\it permutation braid} if its every pair of strands crosses at most once.

A permutation braid is uniquely defined by its permutaion of strands and a number of permutation braids equals to $n!$.

\noindent{\bf Definition.} Let {\it Garside element} be the permutation braid $\Delta$ with $p (i) = n{+}1{-}i$ as a permutation of strands.

For every permutation braid one can find a permutation braid such that their product equals to Garside element.

\noindent{\bf Theorem-definition.} [\ref{processing}] Any braid has a unique representative in {\it left normal form} $\Delta^k\cdot b_1\cdot b_2\cdot\dots\cdot b_m$ where $b_i$ is a permutation braid not equal to $\Delta$ and for $i=1,2,\dots,m{-}1$ and $j=1,2,\dots,n{-}1$ if the $j$th and the $(j{+}1)$th strands of $b_{i+1}$ cross then two strands of $b_i$ which end in the $j$th and the $(j{+}1)$th points also cross.

\noindent{\bf Definition.} [\ref{elrifai}] Let {\it Super Summit Set} of a braid $b$ be the set of braids which are conjugate to $b$ and have maximal power of $\Delta$ and minimal number of permutation braids in their left normal form.

\noindent{\bf Definition.} Consider a braid $b$ in left normal form. Let ${\bf c}(b)=\Delta^k\cdot b_2\cdot b_3\cdot\dots\cdot b_m\cdot(\Delta^k b_1 \Delta^{-k}) =(\Delta^k b_1 \Delta^{-k})^{-1}b(\Delta^k b_1 \Delta^{-k})$ and ${\bf d}(b)=\Delta^k\cdot(\Delta^{-k}b_m\Delta^k)\cdot b_1\cdot b_2\cdot b_{m-1}=b_m b b_m^{-1}.$ We call ${\bf c}(b)$ and ${\bf d}(b)$ the {\it cycling} of $b$ and {\it decycing} of $b$ respectively.

\noindent{\bf Theorem.} [\ref{elrifai}] Let $l$ be the word length of a braid $b$. Then a sequence of at most $l\cdot n^2$ cyclings and decyclings applied to $b$ produces a representative of Super Summit Set.

\noindent{\bf Definition.} [\ref{gebhardt}] Let {\it Ultra Summit Set} ($U_A$) be the subset of Super Summit Set which consists of braids $b$ such that ${\bf c}^d(b){=}b$ for some natural number $d$.

\noindent{\bf Definition.} If ${\bf c}(b)$ in definition of the cycling is already in left normal form then we call $b$ {\it rigid}.

The following definitions refer to theorem 2. Here the main reference is [\ref{bkl}].

For $t>s$ denote the braid $\sigma_{t-1} \sigma_{t-2}\dots\sigma_{s}\sigma_{s+1}^{-1}\sigma_{s+2}^{-1}\dots\sigma_{t-1}^{-1}$ by $(t\ s)$. The set of the braids $(t\ s)$ is called {\it band generators}.

\noindent{\bf Definition.} Let $n_1,n_2,\dots,n_m$ be a decreasing sequence. Denote by $(n_1\ n_2\ \dots\ n_m)$ the braid $(n_1\ n_2)(n_2\ n_3)\dots(n_{m-1}\ n_m)$. This braid is called {\it descending cycle}.

\noindent{\bf Definition.} The braid $\delta:=(n\ n{-}1\dots2\ 1)$ is called a {\it fundamental word}.

\noindent{\bf Definition}. Two descending cycles $(n_1\ n_2\ \dots\ n_p)$ and $(m_1\ m_2\ \dots\ m_q)$ are {\it parallel} if for any $i=1,2,\dots,p{-}1$ and $j=1,2,\dots,q{-}1$ $(n_i-m_j)(n_i-m_{j+1})(n_{i+1}-m_j)(n_{i+1}-m_{j+1})>0.$

\noindent{\bf Definition.} The product of pairwise parallel descending cycles is called a {\it cannonical factor.}

\noindent{\bf Theorem-definition.} Any braid has a unique representation in {\it left normal form} $\delta^k\cdot b_1\cdot b_2\cdot\dots\cdot b_m$ where $b_i$ is a cannonical factor not equal to $\delta$ and for $i=1,2,\dots,m{-}1$ and every generator $(j\ k)$ the braids $b_i\cdot (j\ k)$ and $(j\ k)^{-1}\cdot b_{i+1}$ are not simultaneously cannonical factors.

Definitions if Cycling, Decycling, Super Summit Set and Ultra Summit Set are similar to Artin generators case. Denote Ultra Summit Set in Birman-Ko-Lee presentation by $U_{BKL}$. For definitions of periodic, reducible and pseudo-Anosov braid see [\ref{bestvina},\ref{rigid},\ref{thurston}].

\section{Main result.}

\noindent{\bf Theorem 1.} The braid $\alpha_{n+1}$ on $n{+}1$ strands is rigid and the size of its $U_A$ is at least $2^{[(n-1)/2]}.$

The word generators permutations applied to the braid $\alpha_{n+1}$ produces $2^{n-1}$ braids (proposition 2) which are conjugate to $\alpha_{n+1}$ (proposition 1). We will prove that  $2^{[(n-1)/2]}$ braids among them are rigid. Note that cycling of rigid braid is also rigid. So by the theorem [\ref{elrifai}] mentioned above a rigid braid belongs to its $U_A$. So we obtain $2^{[(n-1)/2]}$ elements of $U_A$.

\noindent{\bf Proposition 1.} A braid produced by generators permutation applied to the word of the braid $\alpha_{n+1}$ is conjugate to $\alpha_{n+1}$.

\noindent{\bf Proof.} Assume that a generator $\sigma_i$ is on the right of the generator $\sigma_1$ in that braid word. Transpose them if $i\neq2$. Do this operation while it is possible. If $\sigma_1$ is on the right end of the word then conjugate the braid by $\sigma_1$ and obtain a braid with $\sigma_1$ being on the left end. Then do mentioned operation until $\sigma_1$ meets $\sigma_2^{-1}$. Then similarly move $\sigma_1$ and $\sigma_2^{-1}$ together on the right until $\sigma_2^{-1}$ meets $\sigma_3$. So at most $n^2$ transpositions and conjugations produce $\alpha_{n+1}$.

\noindent{\bf Proposition 2.} A number of the obtained braids equals to $2^{n-1}$.

\noindent{\bf Proof.} An order of generators of neighbouring numbers in the braid word (which generator $\sigma_i^{\pm1}$ or $\sigma_{i+1}^{\mp1}$ is on the left) determines the braid due to the commutativity relations. We will prove that this order is determined by the braid.

\noindent{\bf Definition.} We correspond to the braid obtained by permuting generators in the word of $\alpha_{n+1}$ a sequence $n_1,m_1,n_2,m_2,\dots,n_r,m_r$ of natural numbers such that $\sigma_i^{\pm1}$ is on the left of $\sigma_{i+1}^{\mp1}$ for $i=1,2,\dots,n_1$, $\sigma_i^{\pm1}$ is on the right of $\sigma_{i+1}^{\mp1}$ for $i=n_1+1,n_1+2,\dots,n_1+m_1$, etc.

The $(n_1{+}m_1{+}\dots{+}n_k{+}m_k{+}1)$th strand ends in the $(n_1{+}m_1{+}\dots{+}n_k{+}m_k{+}n_{k+1}{+}2)$th endpoint and the $(n_1{+}m_1{+}\dots{+}n_k{+}1)$th strand ends in the $(n_1{+}m_1{+}\dots{+}n_k{+}m_k{+}2)$th endpoint. The number of a strand decreases by one if it is between $n_1{+}m_1{+}\dots{+}n_k{+}m_k{+}1$ and $n_1{+}m_1{+}\dots{+}n_k{+}m_k{+}n_{k+1}{+}1$. Other strands increase their number by one. So the order of generators of neighbouring numbers is determined by the sequence $n_1,m_1,\dots,n_r,m_r$ determined by our braid permutation. So we have $2^{n-1}$ braids.

\noindent{\bf Collorary from the proof.} A braid is determined by a sequence $n_1,m_1,n_2,m_2,\dots,n_r,m_r$

\smallskip
\smallskip

According to computations of J.Gonz\'alez-Meneses' program all these braids are rigid (therefore belong to $U_A$ of $\alpha_{n+1}$) for $2\leqslant n\leqslant5$. For simplicity we will consider only such braids among them that if a generator is on the left of both generators of neighbouring numbers then it has an even number and if a generator is on the right of both generators of neighbouring numbers then it has an odd number. In notations of the previous paragraph we require that all $n_i$ are even and all $m_i$ except the last are odd. Also we require that $m_r\neq0$.

\smallskip
\smallskip

\noindent{\bf Proposition 3.} A number of the braids we consider is at least $2^{[\frac{n-1}2]}.$

\noindent{\bf Proof.} Without loss of generality we assume that $n$ is odd. By collorary from proposition 2, a number of the braids we consider equals to a number of decompositions of $n{-}1$ into the sum of $2r$ numbers $n_1,n_2,\dots,n_r,m_1,m_2,\dots,m_r$ for all $r$, where $n_i$ are even and $m_i$ except the last are odd. Assume that $m_r$ is odd. Then the quantity of decompositions we have equals to the quantity of decompositions of $n{-}1{+}r$ into the sum of $2r$ even numbers for all $r$ that equals to the quantity of decompositions of $(n{-}1{+}r)/2$ into the sum of $2r$ numbers for all $r$ that equals to $\sum\limits_{r=1}^{n}C_{(n{-}1{+}r)/2+2r-1}^{2r-1}\geqslant\sum\limits_{r=1}^{n}C_{(n{-}1)/2}^{2r-1}=2^{(n-1)/2-1}.$ For even $m_r$ we obtain the other $2^{(n-1)/2-1}$ braids.

\smallskip

\noindent{\bf Proposition 4.} The braids considered are rigid.

\noindent{\bf Proof.} First we compute a left normal form of the braids considered.

\begin{figure}[h]
\center{\includegraphics[scale=2]{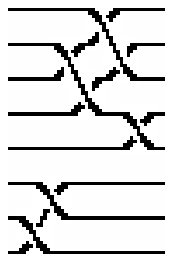} \Huge \  \includegraphics[scale=2]{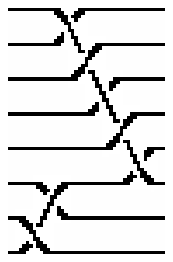} \ \ \ \ \ \  \includegraphics[scale=2]{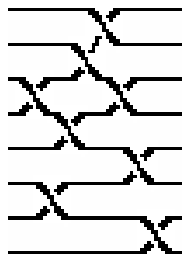}}
\end{figure}

First assume that $r=1$ and $n$ is odd, i.e. we have a braid in the middle of the picture $$\sigma_n\sigma_{n-1}^{-1}\dots\sigma_{n_1+2}^{-1}\sigma_1\sigma_2^{-1}\dots\sigma_{n_1+1}.$$ Note that crossings of the first strand with strands of numbers $3, 5, 7, \dots, n_1{+}1$ are negative. We can get rid of this: multiply our braid from the left by the permutation braid having $(3\ 5\ 7\ \dots\ n_1{+}1\ 1\ 2\ 4\ \dots\ n_1\ n_1{+}3\ n_1{+}5\ \dots\ n_1{+}m_1{+}1\ n_1{+}m_1{+}2\ n_1{+}2\ n_1{+}4 \dots n_1{+}m_1)$ as a permutation of strands. In the obtained product first $n_1/2$ strands cross the other strands from above, which means that we can move them above the other strands to the top, so that negative crossings between first $n_1{+}1$ strands disappear and any two of these strands cross at most once (on the picture the right braid is a product of the braids on the left). Similarly move below to the bottom strands of numbers from $n_1{+}2{+}m_1/2$ to $n$. We obtain a permutation braid having $$(2\ 4\dots n_1\ n_1{+}2\ 1\ 3\ 5\ \dots n_1{+}1\ n_1{+}4\ n_1{+}6 \dots n_1{+}m_1{+}2\ n_1{+}1\ n_1{+}3\ n_1{+}5 \dots n_1{+}m_1{+}1)$$ as a permutation of strands.

Now consider the general case. First reduce our braid by generators which are on the left of both generators of neighbouring numbers: multiply our braid from the left by a product $\sigma_{1+n_1+m_1}\sigma_{1+n_1+m_1+n_2+m_2}\dots\sigma_{1+n_1+m_1+\dots+n_{r-1}+m_{r-1}}$. If $n$ is odd also multiply our braid by $\sigma_n$ from the left. By our assumption these generators have even numbers, therefore they will cancel. Denote that product of generators (including $\sigma_n$ if neccesary) by $c$. After multiplication our braid devides into the product of braids. Each braid among them has all strands standing still except strands whose number is between $1+n_1+m_1+\dots+n_k+m_k$ and $n_1+m_1+\dots+n_{k+1}+m_{k+1}$. Note that the case of such braid was analyzed in the previous paragraph providing us a permutation braid to multiply from the left, which means that we can multiply from the left our braid by a product $d$ of the corresponding permutation braids and obtain a permutation braid $b$. Denote a permutation braid $d\cdot c$ by $a$ and our initial braid by $\beta$. So we have $\beta=a^{-1}b=\Delta^{-1} a^* b$ where $a^*=\Delta a^{-1}.$

Now we prove that $\Delta^{-1} a^* b$ is a left normal form. In the braid $b$ the $i$th and the $(i+1)$th strands cross for $i=n_1/2+1,n_1+m_1/2+2,n_1+m_1+3+n_2/2,n_1+m_1+4+n_2+m_2/2,\dots,n_1+m_1+\dots+n_{r-1}+m_{r-1}+3+n_r/2,n_1+m_1+\dots+n_{r-1}+m_{r-1}+4+n_r+m_r/2.$ But in the braid $a$ for such values of $i$, the $i$th and the $(i+1)$th strands do not cross. Therefore strands ending at the $i$th and the $(i+1)$th endpoints cross in the braid $a^*$. So we have a left normal form.

Now we prove that the braid $\beta$ is rigid. It is sufficient to check that $\Delta^{-1} b (\Delta^{-1} a^* \Delta)$ is a left normal form. Strands ending at the $i$th and the $(i+1)$th endpoints do not cross for odd $i$ in the braid $a$. Note that $a (\Delta^{-1} a^* \Delta)=\Delta$. Therefore the $i$th and the $(i+1)$th strands cross in the braid $(\Delta^{-1} a^* \Delta)$ for odd $i$. But in the braid $b$ strands ending at the $i$th and the $(i+1)$th endpoints crosses for odd $i$. So $\beta$ is rigid.

\bigskip

\noindent{\bf Theorem 2} Ultra Summit Set in Birman-Ko-Lee presentation of the braid $\alpha_n$ contains at least $2^{(n{-}1)/2}$ rigid braids for odd $n$.

%We consider a braid $\sigma_{n-1} \sigma_{n-2}^{-1} \sigma_{n-3}\dots \sigma_1^{(-1)^{n}}=(n\ n{-}1)(n{-}2\ n{-}1)\dots(1\ 2)$ which is conjugate to $\alpha_n$ by $\Delta$.

\noindent{\bf Proposition 1.} A braid $\beta=\delta^{-1}((n{-}1\ 1)(n{-}2\ 2)\dots(\frac{n+1}2\ \frac{n-1}2))^2$ is conjugate to $\alpha_n$.

\noindent{\bf Proof.} A braid $\delta^{-1} \beta \delta=\delta^{-1}( (n\ 2)(n{-}1\ 3)\dots(\frac{n+3}2\ \frac{n+1}2))^2$ is conjugate to the braid $\beta$. So then it is sufficient to prove for this braid the conjugacy to $\alpha_n.$

Using $(n\ 2)(n\ n{-}1\ n{-}2\ \dots 2)=(n\ n{-}1\ n{-}2\ \dots\ 2)(3\ 2)$ we have

\begin{multline*} \bigl((n\ n{-}1\ \dots\ 2) (n\ n{-}1\ \dots\ 4) \dots (n\ n{-}1)\bigr)^{-1} (\delta^{-1}\beta\delta) \times \\ \bigl((n\ n{-}1\ \dots\ 2) (n\ n{-}1\ \dots\ 4) \dots (n\ n{-}1)\bigr)  = \bigl((n\ n{-}1\ \dots\ 2) (n\ n{-}1\ \dots\ 4) \dots (n\ n{-}1)\bigr)^{-1} \times \\ \times \delta^{-1} \times \bigl((n\ 2)(n{-}1\ 3)\dots(\frac{n+3}2\ \frac{n+1}2)\bigr)^2 \times \bigl( (n\ n{-}1\ \dots\ 2) (n\ n{-}1\ \dots\ 4) \dots (n\ n{-}1) \bigr) = \\ = \bigl((n\ n{-}1\ \dots\ 2) (n\ n{-}1\ \dots\ 4) \dots (n\ n{-}1)\bigr)^{-1} \delta^{-1}\times \\ \times \bigl( (n\ n{-}1\ \dots\ 2) (n\ n{-}1\ \dots\ 4) \dots (n\ n{-}1) \bigr) \times \bigl((n\ n{-}1)(n{-}2\ n{-}3)\dots(3\ 2)\bigr)^2 = \\ = \delta^{-1} \bigl((n{-}1\ n{-}2)^{-1} (n{-}1\ n{-}2\ n{-}3\ n{-}4)^{-1}\dots(n{-}1\ n{-}2\ \dots\ 1)^{-1}\bigr)\times \\ \times \bigl( (n\ n{-}1\ \dots\ 2) (n\ n{-}1\ \dots\ 4) \dots (n\ n{-}1) \bigr) \times \bigl((n\ n{-}1)(n{-}2\ n{-}3)\dots(3\ 2)\bigr)^2.\end{multline*}

Then note that $(n{-}1\ n{-}2\ \dots\ 1)^{-1} (n\ n{-}1\ \dots\ 2)= \\ \bigl((2\ 1)^{-1} (n{-}1\ n{-}2\dots\ 2)^{-1}\bigr) \bigl((n{-}1\ n{-}2\ \dots2)(n\ 2)\bigr)=(1\ 2)(n\ 2)$, and $(n\ 2)(n\ n{-}1\ \dots\ 4) = (n\ n{-}1\ \dots\ 4)(4\ 2)$. Using it we continue equality \begin{multline*}\delta^{-1} \bigl( (n{-}1\ n{-}2)^{-1} (n{-}1\ n{-}2\ n{-}3\ n{-}4)^{-1}\dots(n{-}1\ n{-}2\dots1)^{-1}\bigr) \times \\ \times \bigl( (n\ n{-}1\ \dots 2) (n\ n{-}1\ \dots\ 4) \dots (n\ n{-}1) \bigr) \times \\ \times \bigl((n\ n{-}1)(n{-}2\ n{-}3)\dots(3\ 2)\bigr)^2 = \delta^{-1}\bigl( (n{-}2\ n{-}1)(n{-}3\ n{-}2)\dots(1\ 2)\bigr)\times\\\times \bigl((n\ n{-}1)(n{-}1\ n{-}3) (n{-}3\ n{-}5)\dots (4\ 2)\bigr) \times  \\ \times \bigl((n\ n{-}1)(n{-}2\ n{-}3)\dots(3\ 2)\bigr)^2 = \bigl((2\ 3)(4\ 5)\dots(n{-}1\ n)\bigr) \delta^{-1} \times \\ \times (n\ n{-}1\ n{-}3\ n{-}5\dots 2) \bigl((n\ n{-}1)(n{-}2\ n{-}3)\dots(3\ 2)\bigr)^2. \\ \end{multline*}

Conjugating the obtained braid by a braid $(n\ n{-}1)(n{-}2\ n{-}3)\dots(3\ 2)$, we obtain that $\beta$ is conjugate to a braid

\begin{multline*}\bigl(\delta^{-1} (n\ n{-}1\ n{-}3\ n{-}5\ \dots\ 2)\bigr) \bigl((n\ n{-}1)(n{-}2\ n{-}3)\dots(3\ 2)\bigr) = \\ = \Bigl(\bigl((1\ 2)(2\ 3)\dots(n{-}1\ n)\bigr)\bigl((n\ n{-}1)(n{-}1\ n{-}3)(n{-}3\ n{-}5)\dots(4\ 2)\bigr)\Bigr)\times \\ \times\bigl((n\ n{-}1)(n{-}2\ n{-}3)\dots(3\ 2)\bigr)=\\=\bigl((n{-}2\ n{-}1)(n{-}4\ n{-}3)\dots(1\ 2)\bigr)\bigl((n\ n{-}1)(n{-}2\ n{-}3)\dots(3\ 2)\bigr).\end{multline*}

Write the obtained braid in Artin presentation: $\sigma_{n{-}2}^{-1}\sigma_{n-4}^{-1}\dots \sigma_{1}^{-1}\sigma_{n-1}\sigma_{n-3}\dots\sigma_2$. And the same for $\Delta\alpha_n\Delta^{-1}$: $\sigma_{n-1}\sigma_{n-2}^{-1}\sigma_{n-3}\sigma_{n-4}^{-1}\dots\sigma_2\sigma_{1}^{-1}$. Therefore these two braids are conjugate due to Proposition 1 in the proof of Theorem 1.

\smallskip

\noindent{\bf Proposition 2.} Consider a set of braids ${(n{-}1\ 1), (n{-}2\ 2),\dots,(\frac{n+1}2\ \frac{n-1}2)}$. Choose a subset $S$. Denote by $s$ a product of braids in $S$. Then a braid $(\delta^{-1}s\delta)^{-1} \beta (\delta^{-1}s\delta)$ belongs to its $U_{BKL}$ (and therefore to $U_{BKL}$ of $\beta$ and $\alpha_n$).

\noindent{\bf Proof.} Denote by $T$ a complementary to $S$: ${(n{-}1\ 1), (n{-}2\ 2),\dots,(\frac{n+1}2\ \frac{n-1}2)}\setminus S.$ Denote by $t$ a product of braids in $T$. Note that $(\delta^{-1}s\delta)^{-1} \beta (\delta^{-1}s\delta)=\delta^{-1} (ts) (t \delta^{-1}s\delta )$.

\noindent{\bf Step 1.} First we prove that $t \delta^{-1}s\delta $ is a cannonical factor.

Denote $a_i=(n{-}i\ i)$ if $(n{-}i\ i)\notin S$ and $a_i=(n{-}i{+}1\ i{+}1)$ if $(n{-}i\ i)\in S$, for $i=1,2,\dots,\frac{n-1}2$. Consider the braids $(n{-}i\ i)$ belonging to $S$. Then a product of the correspinding $a_i$ equals to $\delta^{-1}s\delta$.

Analyse a product $(t \delta^{-1}s\delta )$. First separate it into the commutating multipliers. Assume that $(\frac{n+1}2\ \frac{n-1}2)\in T$. If $(\frac{n+1}2{+}1\ \frac{n-1}2{-}1)\in T$ then $(\frac{n+1}2\ \frac{n-1}2)$ commutates with the rest of the product, so separate it. If $(\frac{n+1}2{+}1\ \frac{n-1}2{-}1)\in S$ and $(\frac{n+1}2{+}2\ \frac{n-1}2{-}2)\in S$ then $((\frac{n+1}2\ \frac{n-1}2)\cdot(\frac{n+1}2{+}1\ \frac{n-1}2{-}1)$ commutates with the rest of the product, so separate it. Then go further in the same way. Now formalize our observation. We say that $a_i$ and $a_{i+1}$ belong to the same multiplier if $(n{-}i\ i)$ and $(n{-}i-1\ i+1)$ belong to the same subset ($S$ or $T$). Indeed, the product divides into the commutating multipliers in this way.

Now compute each multiplier. It equals to $a_i\cdot a_{i+1}\cdot\dots\cdot a_j$. Its form depends on the belonging of $(n{-}i\ i)$ and $(n{-}j\ j)$ to $S$. So we have four cases. Consider two of them, the other cases are similar. So assume that $(n{-}i\ i)\in T$ and $(n{-}j\ j)\in T$. Then using two formulas ( $(r_1\ r_2\ \dots\ r_p)(r_q\ r)=(r_1\ r_2\dots\ r_q\ r\ r_{q+1}\ \dots\ r_p)$ if $r_q < r< r_{q+1}$ and $(r\ r_q)(r_1\ r_2\ \dots\ r_p)=(r_1\ r_2\dots\ r_{q-1}\ r\ r_q\ \dots\ r_p)$ if $r_{q-1}< r < r_q$) and induction we obtain that the considered multiplier equals to $(n{-}i\ n{-}i{-}2\ \dots\ n{-}j\ j\ j{+}2\ \dots\ i).$ And if $(n{-}i\ i)\in S$ and $(n{-}j\ j)\in T$ then it equals to $(n{-}i{+}1\ n{-}i{-}1\dots n{-}j\ j\ j{+}2\ \dots\ i{+}1).$ In both cases we obtain a descending cycle. Note that the numbers in the cycle form two arithmetic progressions.

Note that the descending cycles in our product are parallel, so the statement of the first step is proved.

\noindent{\bf Step 2.} Now we prove that $\delta^{-1}\cdot(ts)\cdot(t\delta^{-1}s\delta)$ is a left normal form. First we introduce some notations and present some results from [\ref{bkl}].

We say that cannonical factor $c$ is divisible by a generator $(i\ j)$ if $(j\ i)\cdot c$ is also a cannonical factor. Collorary 3.7 [\ref{bkl}]: a cannonical factor is divisible by a generator $(i\ j)$ if and only if one of descending cycles in its decomposition includes $i$ and $j$. By the same collorary, $(j\ i)\cdot c$ is a cannonical factor if and only if $c\cdot(j\ i)$ is a cannonical factor.

For every cannonical factor $a$ $(\delta^{-1}a)$ and $(a\delta^{-1})$ are also cannonical factors ([\ref{bkl}]). Therefore we have: if $ab=\delta$ where $a$ and $b$ are cannonical factors then $a\cdot(i\ j)$ is a cannonical factor if and only if $b$ is divisible by $(i\ j)$. Recall that $ts=((n{-}1\ 1)(n{-}2\ 2)\dots(\frac{n+1}2\ \frac{n-1}2))$ and notice that $$((n{-}1\ 1)(n{-}2\ 2)\dots(\frac{n+1}2\ \frac{n-1}2))\cdot((n\ 1)(n{-}1\ 2)\dots(\frac{n+3}2\ \frac{n-1}2))=\delta.$$

So we can reformulate a statement of step 2: cannonical factor $(t\delta^{-1}s\delta)$ and factor $((n\ 1)(n{-}1\ 2)\dots(\frac{n+3}2\ \frac{n-1}2))$ are not divisible by the same generator simultanously. Indeed, $((n\ 1)(n{-}1\ 2)\dots(\frac{n+3}2\ \frac{n-1}2))$ is divisible only by $(n\ 1),(n{-}1\ 2),\dots,(\frac{n+3}2\ \frac{n-1}2)$. A descending cycle of $(t\delta^{-1}s\delta)$ consists of two arithmetical progressions. Numbers in the first progression are at least $\frac{n+1}2$ and in the second --- at most $\frac{n+1}2$. Therefore if $n{-}i$ and $i{+}1$ belong to some descending cycle then these numbers belong to the different progressions. But recall that a sum of numbers from different progressions is of the same parity with $n$. So the statement of the second step is proved.

\noindent{\bf Step 3.} Now we prove that $(\delta^{-1}ts\delta)\cdot(\delta^{-1}t\delta^{-1}s\delta^2)\cdot\delta^{-1}$ is a right normal form.

We are to prove that if $(\delta^{-1}ts\delta)$ is divisible by $(i\ j)$ then $(i\ j)\cdot t\delta^{-1}s\delta$ is not a cannonical factor. Due to the decomposition into descending cycles of $(\delta^{-1}ts\delta)$ and collorary 3.7 [\ref{bkl}] we have that if $(\delta^{-1}ts\delta)$ is divisible by $(i\ j)$ then $(i\ j)=(n{-}i\ i)$. So denote $r:=(n{-}i\ i)\cdot t\delta^{-1}s\delta$. Assume that $(n{-}i\ i)\in T$. Then the braid $t\delta^{-1}s\delta$ is divisible by $(n{-}i\ i)$ and so $(n{-}i\ i)^{-2} r$ is also a cannonical factor. Therefore $r$ is not a cannonical factor by lemma 3.3 in [\ref{bkl}]. Assume that $(n{-}i\ i)\in S$. Then $((n{-}i\ i)(n{-}i{+}1\ i{+}1))^{-1}r$ is also a cannonical factor. So $r$ is not a cannonical factor by the same lemma.

\noindent{\bf Step 4.} Let $\delta^{-1}\cdot a \cdot b$ and $(\delta^{-1} a \delta)\cdot(\delta^{-1} b \delta)\cdot\delta^{-1}$ be a left normal form and a right normal form respectively. Then a braid $\gamma$ they representing is rigid and so belongs to $U_{BKL}$.

\noindent{\bf Proof.} We use arguments of the second step.  By the definition of cycling ${\bf c}(\gamma)=\delta^{-1} b\cdot (\delta^{-1}a\delta)$. We have $(\delta^{-1} a \delta)\cdot(\delta^{-1} b \delta)\cdot\delta^{-1}$ is a right normal form. Then cannonical factors $(\delta^{-1} a \delta)$ and $b^*$ have not common divisors, where $b^* (\delta^{-1} b \delta)=\delta$. Note that $b\cdot b^* = \delta$. Then $\delta^{-1} b\cdot (\delta^{-1}a\delta)$ is a left normal form and $\gamma$ is rigid. 

Note that cycling of rigid braid is also rigid. So by theorem [\ref{elrifai}] mentioned in the part Definitions a rigid braid belongs to its $U_{BKL}$.

Now finish the prove of proposition 2. In steps 1,2 we found a left normal form of the braid $(\delta^{-1}s\delta)^{-1} \beta (\delta^{-1}s\delta)$. In steps 3,4 we proved that it is rigid and belongs to $U_{BKL}$.

To finish the proof of theorem 2 count a number of braids $(\delta^{-1}s\delta)^{-1} \beta (\delta^{-1}s\delta)$. They are all distinct due to they have distinct left normal forms and their number equals to the number of subsets of $S$, i.e. to $2^{\frac{n-1}2}$. Therefore a size of $U_{BKL}$ is at least $2^{\frac{n-1}2}$.

\bigskip

\noindent{\bf Theorem 3.} The braid $\alpha_n$ is pseudo-Anosov for odd number $n\ge3.$

\noindent{\bf Proof.} Assume that $\alpha_n$ is periodic braid, i.e. $\alpha_n^k$ belongs to center of $B_n$, i.e. it equals to some power of a braid $(n\ n{-}1\ n{-}2\ \dots\ 1)^n$ [Chow, 1948]. But a differ between the quantity of positive generators and negative generators in the word of a braid is invariant with respect to the braid group relations. For the braid $\alpha_n^k$ this differ equals to zero, but for $(n\ n{-}1\ n{-}2\ \dots\ 1)^n$ it is a positive number. Therefore the braid $\alpha_n^k$ is trivial. But this is impossible because the braid group is torsion-free.

Assume that $\alpha_n$ is reducible, i.e. one obtain the braid $\alpha_n$ if one substitute nontrivially the strands of some nontrivial braid by some braids. Denote by $\gamma$ the braid whose strands are substituted. Assume that the $i$th strand is substituted by a braid $\gamma_i$.

Consider a special case $n=3$. Note that $\gamma$ has two strands. Any pair of strands in $\alpha_3^3$ are not linked. Therefore $\gamma^3$ is trivial. Therefore $\gamma_i$ is trivial which means that $\alpha_3^3$ is trivial. Contradiction.

Consider a general case. Let $k,l,m$ be natural numbers, $1\leqslant k< l< m\leqslant n$ and numbers $l{-}k$ and $m{-}l$ are odd. Delete from $\alpha_n^n$ all strands except the $k$th, $l$th and $m$th. Note that we obtain $\alpha_3^3$. Therefore the $k$th, $l$th and $m$th strands belong to either the same $\gamma_i$ or three distinct braids $\gamma_i$.

Assume that the $i$th and $j$th strands belong to the same braid $\gamma_s$ and $i<j$.

Assume that $j{-}i$ is even. In a set of three numbers $i,i{+}1,j$ neighbouring numbers have odd differ. Therefore $(i{+}1)$th strand belongs to $\gamma_s$. Considering sets of three numbers $(i,i{+}1,i{+}2),(i{+}1,i{+}2,i{+}3),\dots,(n{-}2,n{-}1,n)$ and $(i{-}1,i,i{+}1),(i{-}2,i{-}1,i),\dots(1,2,3)\\$ we obtain that all strands belong to $\gamma_s$. That is a contradiction.

Assume that $j{-}i$ is odd. If $i>1$ then consider a set $i{-}1,i,j$ and obtain that $(i{-}1)$th strand belongs to $\gamma_s$. If $j<n$ consider a set $i,j,j{+}1$ and obtain that $(j{+}1)$th strand belongs to $\gamma_s$. Then we obtain similarly that all strands belong to $\gamma_s$. That is a contradiction and so theorem 3 is proved.

\newpage
\newcounter{num}
\setcounter{num}0

\refstepcounter{num}
[\arabic{num}\label{bestvina}] M.Bestvina, M.Handell, Train-tracks for surface homeomorphisms, Topology, vol. 34, 1995, no.1, pp. 1-51

\refstepcounter{num}
[\arabic{num}\label{rigid}] J.Birman, V.Gebhardt, J.Gonz\'alez-Meneses, Conjugacy in Garside groups I: cyclings, powers and rigidity, Groups, Geometry and Dynamics,1,2007,pp.221-279

\refstepcounter{num}
[\arabic{num}\label{periodic}] J.Birman, V.Gebhardt, J.Gonz\'alez-Meneses, Conjugacy in Garside groups III: periodic braids, Journal of Algebra, vol.316, 2007, 2, pp.746-776

\refstepcounter{num}
[\arabic{num}\label{bkl}] J.Birman, K.H.Ko, S.J.Lee, A new approach to the word and conjugacy problems in the braid groups, Advanced Mathematics, vol. 139, 1998, no. 2

\refstepcounter{num}
\noindent[\arabic{num}\label{elrifai}] E.A.Elrifai, H.R.Morton --- Algorithms for positive braids, Quart. J. Math. Oxford (2), 45, 1994, pp. 479-497

\refstepcounter{num}
[\arabic{num}\label{processing}] D.Epstein, J.Cannon, D.Holt, S.Levy, M.Paterson, W.Thurston, Word Processing in Groups, Jones and Barlett Publishers, Boston, MA, 1992

\refstepcounter{num}
[\arabic{num}\label{franco}] N.Franco, J.Gonz\'alez-Meneses, Conjugacy problem for braid groups and Garside groups, Journal of Algebra, 266, 2003, 1, pp. 112-132 

\refstepcounter{num}
[\arabic{num}\label{garside}] F.A.Garside --- The braid group and other groups,Quart. J. Math. Oxford Ser. (2) 20 (1969), 235-254.

\refstepcounter{num}
[\arabic{num}\label{gebhardt}] V.Gebhardt --- A new approach to the conjugacy problem in Garside groups,J. Algebra 292 (1) (2005) 282-302

\refstepcounter{num}
[\arabic{num}\label{thurston}] W.P.Thurston --- On the geometry and dynamics of diffeomorphisms of surfaces, Bull. AMS 19 (1988) no. 2, 417-431

\end{document}